# Alternative to Ritt's Pseudodivision for finding the input-output equations in algebraic structural identifiability analysis

**Nicolette Meshkat[1]\*, Chris Anderson[2], and Joseph J. DiStefano III[3]**

[1,2] UCLA Department of Mathematics, [3] UCLA Department of Computer Science

**Abstract**

Differential algebra approaches to structural identifiability analysis of a dynamic system model in many instances heavily depend upon Ritt's pseudodivision at an early step in analysis. The pseudodivision algorithm is used to find the *characteristic set*, of which a subset, the *input-output equations*, is used for identifiability analysis. A simpler algorithm is proposed for this step, using Gröbner Bases, along with a proof of the method that includes a reduced upper bound on derivative requirements. Efficacy of the new algorithm is illustrated with two biosystem model examples.

## 1. Introduction

*A priori* structural identifiability analysis is concerned with finding one or more sets of solutions for the unknown parameters $\boldsymbol{p}$ of a structured dynamic system model with state and output equations of the form (1.1) from noise-free input-output $\{\boldsymbol{u}(t), \boldsymbol{y}(t)\}$ data:

$$\dot{\boldsymbol{x}}(t, \boldsymbol{p}) = \boldsymbol{f}(\boldsymbol{x}(t, \boldsymbol{p}), \boldsymbol{u}(t), t; \boldsymbol{p}), t \in [t_0, T]$$

$$\boldsymbol{y}(t, \boldsymbol{p}) = \boldsymbol{g}(\boldsymbol{x}(t, \boldsymbol{p}); \boldsymbol{p}) \qquad (1.1)$$

Here $\boldsymbol{x}$ is a n-dimensional state variable, $\boldsymbol{p}$ is a $P$-dimensional parameter vector, $\boldsymbol{u}$ is the r-dimensional input vector, and $\boldsymbol{y}$ is the m-dimensional output vector. In the differential algebra approach, one assumes $\boldsymbol{f}$ and $\boldsymbol{g}$ are rational polynomial functions of their arguments, a reasonable assumption in most applications. The assumption that the output vector is only dependent upon elements of the state variable, and not its derivatives, will be important for the analysis in this work.

Differential algebra approaches have been shown to be quite useful in addressing global as well as local identifiability properties of these models [1-4]; and several differential algebra algorithms have been developed and implemented in available software packages [4-6]. Unfortunately, all are encumbered by computational algebraic complexity or other difficulties, and are limited thus far to relatively low dimensional models [7-9]. To alleviate some of this computational complexity, we describe a procedure that simplifies the task of determining the input-output equations, which is an important early step in preparing the system for identifiability analysis, as considered in [10-14]. The general idea behind the simplified procedure is to use a Gröbner Basis instead of the more cumbersome Ritt's pseudodivision to transform (1.1) into an implicit input-output map involving only the elements and derivatives of $\boldsymbol{y}$ and $\boldsymbol{u}$ along with the parameters $\boldsymbol{p}$, as shown in [15, 16]. We extend the ideas of [15, 16] by finding a stricter bound on the minimum number of derivatives of the equations needed in forming the Gröbner Basis for the multi-output case.

**2. Differential algebra approach to Identifiability in brief**

We now summarize the differential algebra approach to structural identifiability, as well as some differential algebraic concepts. For more details on differential algebra, the reader is referred to [17, 18].

From (1.1), an input-output map is determined in implicit form using a process called Ritt's pseudodivision algorithm [4]. The result of the pseudodivision algorithm is called the *characteristic set* [2]. Since the ideal generated by (1.1) is a prime ideal [19], the characteristic set is a finite "minimal" set of differential polynomials which generate the same differential ideal as that generated by (1.1) [4]. The first m equations of the characteristic set are those independent of the state variables, and form the *input-output relations* [4]:

$$\Psi(y, u, p) = 0 \tag{2.1}$$

The m equations of the input-output relations $\Psi(y, u, p) = 0$ are polynomial equations in $u, \dot{u}, \ddot{u}, \ldots, y, \dot{y}, \ddot{y}, \ldots$, called *differential* polynomials [17], with rational coefficients in the elements of the parameter vector $p$.

For example, a simple first-order model, adapted from [20]:

$$\dot{x} = p_1 x + p_2 u$$

$$y = p_3 x$$

with the chosen ranking $\dot{x} > x > \dot{y} > y > u$ yields an input-output equation, $\Psi(y, u, p) = 0$, of the form:

$$\Psi(y, u, p) = \frac{\dot{y}}{p_3} - \frac{p_1}{p_3} y - p_2 u = 0$$

The characteristic set is in general not unique, but the coefficients $c(p)$ of the input-output equations can be fixed uniquely by normalizing the equations to make them monic, for example, by dividing by $p_2$ [4]:

$$\frac{\dot{y}}{p_2 p_3} - \frac{p_1}{p_2 p_3} y - u = 0$$

Structural identifiability can be determined by testing the injectivity of the coefficients $c(p)$, i.e. the model (1.1) is globally identifiable if and only if $c(p) = c(p^*)$ implies $p = p^*$ for arbitrary $p^*$ [4]. Thus, $\frac{1}{p_2 p_3} = \frac{1}{p_2^* p_3^*}$ and $\frac{p_1}{p_2 p_3} = \frac{p_1^*}{p_2^* p_3^*}$ imply that only $p_1 = p_1^*$ and $p_2 p_3 = p_2^* p_3^*$ can be determined in our example, so the model is unidentifiable.

### 3. Ritt's pseudodivision algorithm

Ritt's pseudodivision is the algorithm that has been more commonly used to find the characteristic set of a prime differential ideal generated by a finite set of differential polynomials [17]. The following procedure follows that in [19].

Let $u_j$ be the leader of a differential polynomial $A_j$, which is the highest ranking derivative of the variables appearing in that polynomial. A polynomial $A_i$ is said to be of *lower rank* than $A_j$ if $u_i < u_j$ or, whenever $u_i = u_j$, the algebraic degree of the leader of $A_i$ is less than the algebraic degree of the leader of $A_j$. A polynomial $A_i$ is *reduced with respect to a polynomial $A_j$* if $A_i$ contains neither the leader of $A_j$ with equal or greater algebraic degree, nor its derivatives. If $A_i$ is not reduced with respect to $A_j$ it can be reduced by using the following pseudodivision algorithm:

(1) If $A_i$ contains the kth derivative $u_j^{(k)}$ of the leader of $A_j$, differentiate $A_j$ k times so its leader becomes $u_j^{(k)}$.
(2) Multiply the polynomial $A_i$ by the coefficient of the highest power of $u_j^{(k)}$ and let $R$ be the remainder of the division of this new polynomial by $A_j^{(k)}$ with respect to the variable $u_j^{(k)}$. Then $R$ is reduced with respect to $A_j^{(k)}$. The polynomial $R$ is called the *pseudoremainder* of the pseudodivision.
(3) The polynomial $A_i$ is replaced by the pseudoremainder $R$ and the process is iterated using $A_j^{(k-1)}$ in place of $A_j^{(k)}$ and so on, until the pseudoremainder is reduced with respect to $A_j$.

This algorithm is applied to a set of differential polynomials, rendering each polynomial reduced with respect to each other, to form an auto-reduced set. The result is a *characteristic set*.

In addition to Ritt's algorithm, a number of other algorithms have been developed to find the full characteristic set, such as the Ritt-Kolchin algorithm [18] and the improved Ritt-Kolchin algorithm [21]. Software implementations of these algorithms can be found in the *diffgrob2* package [5] or the *diffalg* package [6]. However, as noted in [14], "algorithms to find the characteristic set are still under development and the existing software packages do not always work well."

The DAISY program [4] uses Ritt's pseudodivision algorithm to obtain the characteristic set and then the input-output equations, i.e. the first m equations of the characteristic set. While the DAISY program is a useful tool in exploring *global* or *local* identifiability properties of systems, the user may want to obtain the input-output equations for other analyses, e.g. for finding identifiable parameter combinations, as in [9]. Copying the characteristic set from DAISY into a different symbolic algebra package is cumbersome due to syntax differences, especially for large systems.

Alternatively, one could implement Ritt's pseudodivision using any symbolic algebra package, as it requires only low level symbolic operations, e.g. differentiation and polynomial division. While this aspect is good from the standpoint of making few demands on the capabilities of a symbolic software system, it has the negative consequence that the method is time consuming to implement and prone to implementation errors. Since only the input-output equations – and not the full characteristic set – are needed for differential algebra identifiability analysis, a simpler method to obtain just the input-output equations can be quite helpful. We propose an alternative procedure here that utilizes differentiation and Gröbner Bases to ease the implementation difficulties.

**4. Alternative method to find input-output equations**

We obtain the input-output relations by taking a sufficient number of derivatives of the system (1.1), followed by computation of a Gröbner Basis of the new system, similar to the method proposed in [15] and [16]. The main difference between our approach and that proposed in [15] and [16] is that we find a stricter bound on the minimum number of derivatives of the equations needed in forming the Gröbner Basis. Following a minimum number of differentiations of the output equations and corresponding state variable equations, the Buchberger Algorithm is used to eliminate all state variables and derivatives of state variables.

In general, for elimination to work, the number of equations must be *strictly greater* than the number of unknowns, as discussed in [15]. Since the output equation is of the form $\boldsymbol{y}(t,\boldsymbol{p}) = \boldsymbol{g}(\boldsymbol{x}(t,\boldsymbol{p});\boldsymbol{p})$, i.e. always in terms of $\boldsymbol{x}$ and not derivatives of $\boldsymbol{x}$, then the first step is to take the derivative of the output equations, to help eliminate the first derivative of $\boldsymbol{x}$ from the state variable equations. If this additional equation is not enough to eliminate the state variables, then the second derivative of the output equations is needed. This, however, introduces the second derivative of $\boldsymbol{x}$, and thus differentiation of the corresponding state variable equations is needed. Differentiation of the output equations and corresponding state variable equations is continued until the number of equations is greater than the number of unknowns. The procedure is described by the following steps:

**Step 1**: Differentiate the output equations to obtain $\dot{\boldsymbol{y}}$ and adjoin these equations to the system.

**Step 2**: Differentiate the output equations again to obtain $\ddot{\boldsymbol{y}}$ and differentiate the corresponding state variable equations to obtain equations involving $\ddot{\boldsymbol{x}}$. Adjoin these equations to the system.

…

**Step k**: Differentiate the output equations again to obtain $\boldsymbol{y}^{(k)}$ and differentiate the corresponding state variable equations to obtain equations involving $\boldsymbol{x}^{(k)}$. Adjoin these equations to the system.

For a single output system, the above steps yield a procedure similar to that in [15] and [16]. However, the method described in [15] and [16] does not formally treat multi-output models, and the procedure we present here does.

We now show that, for a system of n state variables and m output equations, one need only take n-(m-1) steps to obtain a sufficient number of additional equations to eliminate the state variable terms in a Gröbner Basis. We begin with examples illustrating the procedure and complications that arise for more than one output equation. Following the examples, we give a proof that guarantees that the procedure will be successful in constructing the required number of input-output equations.

## 5. Motivating examples

**Example 1.** We demonstrate the process first on a Linear 3-Compartment Model, with one output equation:

$$\dot{x}_1 = k_{13}x_3 + k_{12}x_2 - (k_{21} + k_{31})x_1 + u$$

$$\dot{x}_2 = k_{21}x_1 - (k_{12} + k_{02})x_2$$

$$\dot{x}_3 = k_{31}x_1 - (k_{13} + k_{03})x_3$$

$$y = \frac{x_1}{v}$$

$k_{12}, k_{21}, k_{13}, k_{31}, k_{02}, k_{03}, v$ are the unknown parameters.

**Step 1:** Differentiate the output equation and add this to the system:

$$v\dot{y} = \dot{x}_1$$

We now have 5 equations in 6 unknown state variables and their derivatives, $\{\dot{x}_1, \dot{x}_2, \dot{x}_3, x_1, x_2, x_3\}$, not enough equations to eliminate the unknowns. An additional differentiation of the output equation is thus needed. This introduces the unknown $\ddot{x}_1$, so we must differentiate the $\dot{x}_1$ equation as well, as in Step 2:

**Step 2:**

$$v\ddot{y} = \ddot{x}_1$$

$$\ddot{x}_1 = k_{13}\dot{x}_3 + k_{12}\dot{x}_2 - (k_{21} + k_{31})\dot{x}_1 + \dot{u}$$

This gives 7 equations in 7 unknown state variables and their derivatives, $\{\ddot{x}_1, \dot{x}_1, \dot{x}_2, \dot{x}_3, x_1, x_2, x_3\}$, still not enough equations to eliminate the unknowns since, as discussed in [15], the number of equations must be *strictly greater* than the number of unknowns. An additional differentiation of the output equation introduces the unknown $\dddot{x}_1$. We therefore differentiate the $\ddot{x}_1$ equation. This introduces the unknowns $\ddot{x}_2, \ddot{x}_3$, requiring differentiation of the $\dot{x}_2$ and $\dot{x}_3$ equations, as in Step 3:

**Step 3:**

$$v\dddot{y} = \dddot{x}_1$$

$$\dddot{x}_1 = k_{13}\ddot{x}_3 + k_{12}\ddot{x}_2 - (k_{21} + k_{31})\ddot{x}_1 + \ddot{u}$$

$$\ddot{x}_2 = k_{21}\dot{x}_1 - (k_{12} + k_{02})\dot{x}_2$$

$$\ddot{x}_3 = k_{31}\dot{x}_1 - (k_{13} + k_{03})\dot{x}_3$$

We now have 11 equations in 10 unknown state variables and their derivatives, $\{\dddot{x}_1, \ddot{x}_1, \ddot{x}_2, \ddot{x}_3, \dot{x}_1, \dot{x}_2, \dot{x}_3, x_1, x_2, x_3\}$, enough to eliminate all unknowns. We thus find the Gröbner Basis of the new system, using Mathematica, according to the ranking $\{\dddot{x}_1, \ddot{x}_1, \ddot{x}_2, \ddot{x}_3, \dot{x}_1, \dot{x}_2, \dot{x}_3, x_1, x_2, x_3, \dddot{y}, \ddot{y}, \dot{y}, y, \ddot{u}, \dot{u}, u\}$. We get an equation only in $\dddot{y}, \ddot{y}, \dot{y}, y, \ddot{u}, \dot{u}, u$:

$$(k_{02}k_{03} + k_{03}k_{12} + k_{02}k_{13} + k_{12}k_{13})u + (k_{02} + k_{03} + k_{12} + k_{13})\dot{u} + \ddot{u}$$
$$- (k_{02}k_{03}k_{21} + k_{02}k_{13}k_{21} + k_{02}k_{03}k_{31} + k_{03}k_{12}k_{31})vy$$
$$- (k_{02}k_{03} + k_{03}k_{12} + k_{02}k_{13} + k_{12}k_{13} + k_{02}k_{21} + k_{03}k_{21} + k_{13}k_{21} + k_{02}k_{31}$$
$$+ k_{03}k_{31} + k_{12}k_{31})v\dot{y} - (k_{02} + k_{03} + k_{12} + k_{13} + k_{21} + k_{31})v\ddot{y} - v\dddot{y} = 0$$

This is the input-output equation. The identifiability analysis of this example is demonstrated in [8].

Notice that, for only one output equation, only one input-output equation results. In this case, we needed the number of equations to be strictly greater than the number of unknowns. In our next example, the situation gets more tricky when there are multiple output equations.

**Example 2.** Consider the HIV kinetics model from [22], which has two output equations:

$$\dot{x}_1 = -\beta x_1 x_4 - d x_1 + s$$

$$\dot{x}_2 = \beta q_1 x_1 x_4 - k_1 x_2 - \mu_1 x_2$$

$$\dot{x}_3 = \beta q_2 x_1 x_4 + k_1 x_2 - \mu_2 x_3$$

$$\dot{x}_4 = -c x_4 + k_2 x_3$$

$$y_1 = x_1$$

$$y_2 = x_4$$

$\beta, d, s, q_1, k_1, \mu_1, q_2, k_2, \mu_2, c$ are the unknown parameters.

**Step 1:** Differentiate the output equations and add these to the system:

$$\dot{y}_1 = \dot{x}_1$$

$$\dot{y}_2 = \dot{x}_4$$

We have 8 equations in 8 unknown state variables and their derivatives, $\{\dot{x}_1, \dot{x}_2, \dot{x}_3, \dot{x}_4, x_1, x_2, x_3, x_4\}$, not enough to eliminate the unknowns.

**Step 2:** The second step results in:

$$\ddot{y}_1 = \ddot{x}_1$$

$$\ddot{y}_2 = \ddot{x}_4$$

$$\ddot{x}_1 = -\beta \dot{x}_1 x_4 - \beta \dot{x}_4 x_1 - d \dot{x}_1$$

$$\ddot{x}_4 = -c \dot{x}_4 + k_2 \dot{x}_3$$

This gives 12 equations and 10 unknowns, $\{\ddot{x}_1, \ddot{x}_4, \dot{x}_1, \dot{x}_2, \dot{x}_3, \dot{x}_4, x_1, x_2, x_3, x_4\}$. However, a Gröbner Basis of this system, obtained using Mathematica or a similar program, only gives 1 input-output equation, not the required 2 input-output equations! An input-output equation involving $\dot{y}_1$ is easily obtained, i.e. $\dot{y}_1 + \beta y_1 y_2 + dy_1 - s = 0$, but the input-output equation involving $\dot{y}_2$ cannot yet be obtained. One more round, i.e. Step 3, is needed.

**Step 3**:

$$\dddot{y}_1 = \dddot{x}_1$$

$$\dddot{y}_2 = \dddot{x}_4$$

$$\dddot{x}_1 = -\beta \ddot{x}_1 x_4 - 2\beta \dot{x}_1 \dot{x}_4 - \beta \ddot{x}_4 x_1 - d\ddot{x}_1$$

$$\dddot{x}_4 = -c\ddot{x}_4 + k_2 \ddot{x}_3$$

$$\ddot{x}_3 = \beta q_2 \dot{x}_1 x_4 + \beta q_2 \dot{x}_4 x_1 + k_1 \dot{x}_2 - \mu_2 \dot{x}_3$$

We now have 17 equations and 13 unknowns, $\{\dddot{x}_1, \dddot{x}_4, \ddot{x}_1, \ddot{x}_3, \ddot{x}_4, \dot{x}_1, \dot{x}_2, \dot{x}_3, \dot{x}_4, x_1, x_2, x_3, x_4\}$. We can now find the Gröbner Basis and the required two input-output equations:

$$\dot{y}_1 + \beta y_1 y_2 + dy_1 - s = 0$$

$$\dddot{y}_2 + (c + k_1 + \mu_1 + \mu_2)\ddot{y}_2 - \beta q_2 k_2 \dot{y}_2 y_1 + (ck_1 + c\mu_1 + c\mu_2 + k_1 \mu_2 + \mu_1 \mu_2)\dot{y}_2 + \beta^2 q_2 k_2 y_1 y_2^2 \\ + \beta k_2 (dq_2 - k_1 q_1 - k_1 q_2 - \mu_1 q_2) y_1 y_2 + (-\beta q_2 k_2 s + c k_1 \mu_2 + c \mu_1 \mu_2) y_2 = 0$$

The identifiability analysis of this example is demonstrated in [8].

In the procedure that we have described, the critical fact that is necessary for success is that we must perform up to Step n-(m-1), where n is the number of state variables and m is the number of outputs. In our examples, we showed that for the linear 3-Compartment model with 1 output equation, it was necessary to perform 3 steps, and in the 4-dimensional HIV model with 2 output equations, it was necessary to perform 3 steps as well. In the case of the single output equation, it was shown in [15] and [16] that derivatives up to the nth degree must be taken, in other words Step n will guarantee success of the procedure for the single output equation case. For multiple output equations, we prove below that performing up to Step n-(m-1) is sufficient to generate enough equations for n state variable equations and m output equations.

**6. Number of derivatives of output equations needed**

It was shown in [15] and [16] that for a single output system, an input-output equation of the same (or possibly lower) differential order as the number of state variables can always be obtained. In other words, for a single output system with 4 state variables, an input-output equation of differential order 4 or lower can always be obtained.

This result can be generalized to systems with more than one output equation, simply by examining the state variable equations with each output equation taken individually. For example, in the HIV model above, we can examine the 4 state variable equations along with the second output equation. According to [15] and [16], we should take derivatives up to the fourth order in the output variable in order to obtain an input-output equation. This is likewise true for the first output equation taken with the state space system. So, we can obtain the required two input-output equations in this way, where each input-output equation involves either the first or second output variable, but not both. Thus, an *upper bound* for the number of derivatives of the output equations needed is certainly the number of state variables, n. However, including an additional output equation in our single output system adds more information, thus one is led to conjecture that fewer than n derivatives of the output equations are needed. In the case of the HIV model, only 3 derivatives of the output variable were needed to obtain the two required input-output equations. This is no accident, as we will prove that only 4-(2-1)=4-1=3 derivatives of the output equations are needed in this case. We will prove below that for two output equations, n-1 derivatives of the output equations are needed to obtain the required two input-output equations. We will then, by induction, show that for m output equations, n-(m-1) derivatives of the output equations are needed to obtain the required m input-output equations.

Let our system be of the following form:

$$\dot{x}_1 = f_1(x_1, x_2, \ldots, x_n, u) \tag{6.1}$$

$$\dot{x}_2 = f_2(x_1, x_2, \ldots, x_n, u)$$

$$\ldots$$

$$\dot{x}_n = f_n(x_1, x_2, \ldots, x_n, u)$$

$$y_1 = g_1(x_1, x_2, \ldots, x_n)$$

$$y_2 = g_2(x_1, x_2, \ldots, x_n)$$

We first examine the total number of equations and unknowns at each step. We show that this number varies depending on the particular state variables in the arguments of the functions $f_1, f_2, \ldots, f_n, g_1$, and thus we should instead examine the difference between the number of unknowns and equations at each step.

As mentioned above, we ignore the second output equation. Thus, after Step 1 we have the following:

$$\dot{x}_1 = f_1(x_1, x_2, \ldots, x_n, u) \tag{6.2}$$

$$\dot{x}_2 = f_2(x_1, x_2, \ldots, x_n, u)$$

$$\ldots$$

$$\dot{x}_n = f_n(x_1, x_2, \ldots, x_n, u)$$

$$y_1 = g_1(x_1, x_2, \ldots, x_n)$$

$$\dot{y}_1 = \sum_{i=1}^{n} \frac{\partial g_1}{\partial x_i} \dot{x}_i$$

Thus there are n+2 equations and 2n unknowns, so for n>1 the number of unknowns is greater than or equal to the number of equations.

In Step 2, we take an additional derivative of the output equation, which causes the second derivative of the state variable to appear. Thus, we take the corresponding derivatives of the state variable equations and add these to the system:

$$\ddot{y}_1 = \sum_{i=1}^{n} \left[ \frac{\partial^2 g_i}{\partial x_i \partial t} \dot{x}_i + \frac{\partial g_i}{\partial x_i} \ddot{x}_i \right]$$

$$\ddot{x}_1 = \sum_{i=1}^{n} \frac{\partial f_1}{\partial x_i} \dot{x}_i + \frac{\partial f_1}{\partial u} \dot{u}$$

$$\ddot{x}_2 = \sum_{i=1}^{n} \frac{\partial f_2}{\partial x_i} \dot{x}_i + \frac{\partial f_2}{\partial u} \dot{u}$$

...

$$\ddot{x}_n = \sum_{i=1}^{n} \frac{\partial f_n}{\partial x_i} \dot{x}_i + \frac{\partial f_n}{\partial u} \dot{u}$$

In doing this, we add n+1 equations and n unknowns (if $y_1$ is a function of all n state variables), so that the total is now 2n+3 equations and 3n unknowns. For n>2, the number of unknowns is greater than or equal to the number of equations.

In Step 3, we take an additional derivative of the output equation, which causes the third derivative of the state variable to appear. Thus, we take the corresponding derivatives of the state variable equations and add these to the system:

$$\dddot{y}_1 = \sum_{i=1}^{n} \left[ \frac{\partial^3 g_i}{\partial x_i \partial^2 t} \dot{x}_i + 2 \frac{\partial^2 g_i}{\partial x_i \partial t} \ddot{x}_i + \frac{\partial g_i}{\partial x_i} \dddot{x}_i \right]$$

$$\dddot{x}_1 = \frac{\partial^2 f_1}{\partial u \partial t} \dot{u} + \frac{\partial f_1}{\partial u} \ddot{u} + \sum_{i=1}^{n} \left[ \frac{\partial^2 f_1}{\partial x_i \partial t} \dot{x}_i + \frac{\partial f_1}{\partial x_i} \ddot{x}_i \right]$$

$$\dddot{x}_2 = \frac{\partial^2 f_2}{\partial u \partial t} \dot{u} + \frac{\partial f_2}{\partial u} \ddot{u} + \sum_{i=1}^{n} \left[ \frac{\partial^2 f_2}{\partial x_i \partial t} \dot{x}_i + \frac{\partial f_2}{\partial x_i} \ddot{x}_i \right]$$

...

$$\dddot{x}_n = \frac{\partial^2 f_n}{\partial u \partial t}\dot{u} + \frac{\partial f_n}{\partial u}\ddot{u} + \sum_{i=1}^{n}\left[\frac{\partial^2 f_n}{\partial x_i \partial t}\dot{x}_i + \frac{\partial f_n}{\partial x_i}\ddot{x}_i\right]$$

In doing this, we add n+1 equations and n unknowns, so that the total is now 3n+4 equations and 4n unknowns. For n>3, the number of unknowns is greater than or equal to the number of equations. For example, for n=4 the number of equations and unknowns is equal at this step.

Notice that, in Step 2, if $y_1$ is not a function of all the state variables, for instance, if $y_1 = g_1(x_2, \dots, x_n)$, then only the derivatives

$$\ddot{y}_1 = \sum_{i=2}^{n}\left[\frac{\partial^2 g_i}{\partial x_i \partial t}\dot{x}_i + \frac{\partial g_i}{\partial x_i}\ddot{x}_i\right]$$

$$\ddot{x}_2 = \sum_{i=1}^{n}\frac{\partial f_2}{\partial x_i}\dot{x}_i + \frac{\partial f_2}{\partial u}\dot{u}$$

$$\dots$$

$$\ddot{x}_n = \sum_{i=1}^{n}\frac{\partial f_n}{\partial x_i}\dot{x}_i + \frac{\partial f_n}{\partial u}\dot{u}$$

would be taken at this step. At Step 1 there are always n+2 equations and 2n unknowns, thus there are now (n+2)+n=2n+2 equations and 2n+(n-1)=3n-1 unknowns. However, this means at Step 3, we must have:

$$\dddot{y}_1 = \sum_{i=2}^{n}\left[\frac{\partial^3 g_i}{\partial x_i \partial^2 t}\dot{x}_i + 2\frac{\partial^2 g_i}{\partial x_i \partial t}\ddot{x}_i + \frac{\partial g_i}{\partial x_i}\dddot{x}_i\right]$$

$$\dddot{x}_2 = \frac{\partial^2 f_2}{\partial u \partial t}\dot{u} + \frac{\partial f_2}{\partial u}\ddot{u} + \sum_{i=1}^{n}\left[\frac{\partial^2 f_2}{\partial x_i \partial t}\dot{x}_i + \frac{\partial f_2}{\partial x_i}\ddot{x}_i\right]$$

$$\dots$$

$$\dddot{x}_n = \frac{\partial^2 f_n}{\partial u \partial t}\dot{u} + \frac{\partial f_n}{\partial u}\ddot{u} + \sum_{i=1}^{n}\left[\frac{\partial^2 f_n}{\partial x_i \partial t}\dot{x}_i + \frac{\partial f_n}{\partial x_i}\ddot{x}_i\right]$$

$$\ddot{x}_1 = \sum_{i=1}^{n}\frac{\partial f_1}{\partial x_i}\dot{x}_i + \frac{\partial f_1}{\partial u}\dot{u}$$

More precisely, since $\ddot{x}_1$ now appears, this means we must now include the $\ddot{x}_1$ equation. Thus, there are now (2n+2)+(n+1)=3n+3 equations and (3n-1)+n=4n-1 unknowns.

Notice, in either case, the conserved quantity is that there is one more equation than unknown added to the system at each step. Thus, instead of examining the total number of equations and unknowns at each step, which we have shown varies with the particular state variables in the arguments of the functions $f_1, f_2, \ldots, f_n, g_1$, we will examine the difference between the number of unknowns and equations at each step.

## 7. Proof for new bound

We now present a collection of lemmas and theorems that demonstrate a stricter bound for the number of derivatives of the output equations.

**Lemma 4.1:** Let the system consist of n state variables and 1 output equation. At each step of the algorithm, one more equation than unknown is added to the system.

Proof: Assume the output equation involves a subset of $j$ state variables from the set of state variables, $\{x_1, \ldots, x_n\}$, where $1 \leq j \leq n$. Call this subset of state variables S.

At Step 1: The derivative of the output equation increases the number of equations by 1. No new unknowns are introduced. The net result is an addition of one more equation than unknown.

At Step 2: The second derivative of the output equation increases the number of equations by 1, but $j$ unknowns are introduced. These are the second derivatives of the state variables in S. Adjoin to the system the expressions for the second derivatives of these state variables, in terms of first derivatives of the state variables. The net result is an addition of one more equation than unknown.

At Step k: The kth derivative of the output equation furnishes an additional equation. A collection of $j$ new unknowns are introduced, corresponding to the kth time derivatives of the $j$ state variables occurring in the output equation. Adjoin to this system the expressions for the kth time derivative for these $j$ state variables, in terms of the (k-1)st time derivatives of the state variables. As a result, a set of $i$ new state variables are introduced, which are (k-1)st time derivatives of the state variables, where $0 \leq i \leq n - j$. Adjoin to this system the expressions for the (k-1)st time derivatives for these $i$ state variables, in terms of the (k-2)nd time derivatives of the state variables. Continue this process of adjoining expressions for the state variables in terms of lower order derivatives until no new derivatives of state variables are introduced. As described above, there are at most n+1 equations added at this step, and at most n unknowns. To show that this step results in precisely one more equation than unknown, we use induction.

Assume this is true for Step k-1, i.e. there is one more equation than unknown added at this step. Now take a derivative of all of the equations adjoined to the system in Step k-1. We achieve the same system as before, but each term is of one higher differential order. To achieve the result of the application of Step k, then if any new time derivatives of state variables appear on the right hand side of this system, then the corresponding equations for the state variable derivatives are adjoined to the system. This process will result in an identical system to the system from Step k since the output equation will be of differential order k and any new unknowns added at this step will have the corresponding equations in terms of lower order derivatives. More precisely, for the $r$ new unknowns added to the system as a result of the differentiation, we add the corresponding $r$ equations, where $0 \leq r \leq n - j$. Since we assume at Step k-1 that there is one more equation than unknown, then the net result is one more equation than unknown after this step. Thus, the lemma is proven by induction. ∎

**Lemma 4.2:** Let the system consist of n state variables and 1 output equation. There are always an equal number of equations and unknowns when there are n-1 derivatives of the output equation, i.e. at Step n-1.

Proof: Initially, there are n+1 equations and 2n unknowns. At Step 1, there are always n+2 equations and 2n unknowns, as in (6.2).

Thus there are n-2 more unknowns than equations. In each step thereafter, we add one more equation than unknown, so at the second step, there are n-3 more unknowns than equations. Thus at the n-1$^{st}$ step, there are n-(n-1+1)=0 more unknowns than equations. Thus, the Lemma is proven. ∎

**Corollary to Lemma 4.2:** At the nth step, there is one more equation than unknown.

Proof: This follows from Lemma 4.1 and 4.2. ∎

Remark: This agrees with the result in [15] and [16], i.e., for a single output system, an input-output equation of the same (or possibly lower) differential order as the number of state variables can always be obtained. In other words, since there is one more equation than unknown at Step n, an input-output equation of differential order n, or lower, can be obtained. We will now extend the work of [15] and [16] to formally treat the multiple output case.

**Theorem 4.1:** For a two output system, two input-output equations of differential order one less than the number of state variables (or possibly lower) can always be obtained, i.e. n-1 or lower.

Proof: We showed in Lemma 4.2 that there are an equal number of equations and unknowns at the n-1$^{st}$ derivative for a system with a single output equation. Thus, our state variable system with the first output equation has an equal number of equations and unknowns at the n-1$^{st}$ step, likewise for the system taken with the second output equation. Let System A be the system with the first output equation $y_1 = g_1(x_1, x_2, \ldots, x_n)$ after Step n-1 is performed. Let System B be the system with the second output equation $y_2 = g_2(x_1, x_2, \ldots, x_n)$ after Step n-1 is performed, i.e. System B is the same exact system, except $y_1$ is replaced by $y_2$, $g_1$ is replaced by $g_2$, etc. Thus both System A and System B have an equal number of equations and unknowns.

Now, we will add the equation $y_2 = g_2(x_1, x_2, \ldots, x_n)$ to System A, and we call the new system A'. This adds one new equation, but no new unknowns since $x_1, x_2, \ldots, x_n$ are assumed to already be included in System A. We can do the same for System B, i.e. System B' is System B plus the output equation $y_1 = g_1(x_1, x_2, \ldots, x_n)$. So, for either system, adding an additional output equation is enough to result in more equations than unknowns for each of the systems, because this adds one more equation, but no new unknowns. Thus, the state variables can be eliminated in both Systems A' and B', and the desired two input-output equations of lowest differential order can be obtained. Combining these two systems,

A' and B', will result in an equivalent system to the one obtained by performing Steps 1, 2, …, n-1 on the original system with both output equations. The resulting System A'∪B' will thus result in enough equations for the Buchberger Algorithm to eliminate all the state variables and provide two input-output equations. ∎

Note that in the proof above, we add the output equation $y_2 = g_2(x_1, x_2, \ldots, x_n)$ to System A to obtain System A'. However, it is also acceptable to add $\dot{y}_2 = \frac{\partial g_2}{\partial x_1}\dot{x}_1 + \frac{\partial g_2}{\partial x_2}\dot{x}_2 + \cdots + \frac{\partial g_2}{\partial x_n}\dot{x}_n$ instead of $y_2 = g_2(x_1, x_2, \ldots, x_n)$, since this also adds one more equation but no new unknowns. We note, however, that the resulting input-output equation for this System A' could be a much more complicated equation in terms of $\dot{y}_2$ as opposed to $y_2$. Thus, when we form the System A'∪B' and find the Gröbner Basis, the output equation $y_2 = g_2(x_1, x_2, \ldots, x_n)$ will be used to form the input-output equation associated with System A'.

Since the equations $y_2 = g_2(x_1, x_2, \ldots, x_n)$ and $\dot{y}_2 = \frac{\partial g_2}{\partial x_1}\dot{x}_1 + \frac{\partial g_2}{\partial x_2}\dot{x}_2 + \cdots + \frac{\partial g_2}{\partial x_n}\dot{x}_n$ add two new equations and no new unknowns, one may ask why we perform up to Step n-1 and add one output equation as opposed to Step n-2 and add both equations to form System A'. In the HIV model in Example 2, it was shown that performing up to Step 2 only resulted in one input-output equation. In other words, if we only perform up to Step n-2, then the Systems A' and B' may result in the same input-output equation. This is because each system will result in an input-output equation that is up to differential order n-2 in $y_1$ or $y_2$, but at most differential order 1 in the other output variable. As shown in the HIV model, it is possible to have an input-output equation that is at most differential order 1 in both $y_1$ and $y_2$. Thus, to guarantee two distinct input-output equations in the System A'∪B', we must perform steps up to n-1 since each system (A' or B') will result in an input-output equation that is up to differential order n-1 in $y_1$ or $y_2$, but at most differential order 0 in the other. Thus, the only way to get the same input-output equation for both systems is if the input-output equation is only in terms of $y_1$ and $y_2$ (without derivatives), which is not possible for our ODE system because this means the state variables are algebraically dependent.

We now generalize our result to m output equations:

**Theorem 4.2**: Let the system consist of n state variables and m output equations. Then exactly m input-output equations of differential order m-1 less than the number of state variables (or possibly lower) can always be obtained, i.e. only steps up to n-(m-1) must be performed.

Proof: We have shown that for m=1, we need n derivatives of the output equations to eliminate state variables, and for m=2, we need n-1 derivatives of the output equations to eliminate state variables. We prove this theorem following the logic of the proof in Theorem 4.1.

We first examine the system with only one output equation. We showed in Lemma 4.2 that at the n-1[st] step, there are an equal number of equations and unknowns. This means at the n-2[nd] step, there is one less equation than unknowns, and at the n-3[rd] step, there are two less equations than unknowns, etc. So at the n-(m-1) step, there are m-2 less equations than unknowns. Thus, if we add an additional m-1 output equations, i.e. the other m-1 output equations in our system, then there will be -(m-2)+m-1=1 more equation than unknowns. The same is true for the system taken with any other output equation, thus the theorem is proven. ∎

Thus we have found a minimal number of derivatives to take of the output and state variable equations in order to get the required m input-output equations.

## 8. Additional Example

We now demonstrate our algorithm on another multi-output system. The following model is derived from Evans and Chappell [23] and describes the pharmacokinetics of bromosulphthalein. We find the input-output equations for a similar problem, with a more general input $u(t)$ rather than an initial condition:

$$\dot{x}_1 = -a_{31}x_1 + a_{13}x_3 + u$$

$$\dot{x}_2 = -a_{42}x_2 + a_{24}x_4$$

$$\dot{x}_3 = a_{31}x_1 - (a_{03} + a_{13} + a_{43})x_3$$

$$\dot{x}_4 = a_{42}x_2 + a_{43}x_3 - (a_{04} + a_{24})x_4$$

$$y_1 = x_1$$

$$y_2 = x_2$$

$a_{03}, a_{04}, a_{13}, a_{24}, a_{31}, a_{42}, a_{43}$ are the unknown parameters.

In this model, n=4 and m=2. Thus, to obtain two input-output equations, steps up to 4-(2-1)=4-1=3 must be performed. In other words, we find a Gröbner Basis of the following system:

$$\dot{x}_1 = -a_{31}x_1 + a_{13}x_3 + u$$

$$\dot{x}_2 = -a_{42}x_2 + a_{24}x_4$$

$$\dot{x}_3 = a_{31}x_1 - (a_{03} + a_{13} + a_{43})x_3$$

$$\dot{x}_4 = a_{42}x_2 + a_{43}x_3 - (a_{04} + a_{24})x_4$$

$$y_1 = x_1$$

$$y_2 = x_2$$

$$\dot{y}_1 = \dot{x}_1$$

$$\dot{y}_2 = \dot{x}_2$$

$$\ddot{y}_1 = \ddot{x}_1$$

$$\ddot{y}_2 = \ddot{x}_2$$

$$\ddot{x}_1 = -a_{31}\dot{x}_1 + a_{13}\dot{x}_3 + \dot{u}$$

$$\ddot{x}_2 = -a_{42}\dot{x}_2 + a_{24}\dot{x}_4$$

$$\ddot{y}_1 = \ddot{x}_1$$

$$\ddot{y}_2 = \ddot{x}_2$$

$$\dddot{x}_1 = -a_{31}\ddot{x}_1 + a_{13}\ddot{x}_3 + \ddot{u}$$

$$\dddot{x}_2 = -a_{42}\ddot{x}_2 + a_{24}\ddot{x}_4$$

$$\ddot{x}_3 = a_{31}\dot{x}_1 - (a_{43} + a_{13} + a_{03})\dot{x}_3$$

$$\ddot{x}_4 = a_{42}\dot{x}_2 + a_{43}\dot{x}_3 - (a_{04} + a_{24})\dot{x}_4$$

Thus there are 18 equations and 14 unknowns. We use the following ranking of parameters: $\{\dddot{x}_1, \dddot{x}_2\ \ddot{x}_1, \ddot{x}_2, \ddot{x}_3, \ddot{x}_4, \dot{x}_1, \dot{x}_2, \dot{x}_3, \dot{x}_4, x_1, x_2, x_3, x_4, \ddot{y}_1, \ddot{y}_2, \dot{y}_1, \dot{y}_2, y_1, y_2, \ddot{u}, \dot{u}, u\}$,

which results in the following input-output equations:

$$a_{13}\ddot{y}_2 + (a_{42}a_{13} + a_{13}a_{04} + a_{13}a_{24})\dot{y}_2 + a_{13}a_{42}a_{04}y_2 - a_{43}a_{24}\dot{y}_1 - a_{43}a_{31}a_{24}y_1 + a_{24}a_{43}u = 0$$

$$\ddot{y}_1 + (a_{31} + a_{03} + a_{13} + a_{43})\dot{y}_1 + (a_{31}a_{03} + a_{31}a_{43})y_1 - \dot{u} + (a_{03} + a_{13} + a_{43})u = 0$$

The first input-output equation can be made monic by dividing by $a_{24}a_{43}$, and then identifiability analysis can be performed, as in [8].

**9. Conclusion**

We have presented and exemplified a new algorithm for finding the input-output equations of any nonlinear ODE system of the form (1.1) by taking derivatives of the system (1.1) and then finding a Gröbner Basis. We also found a stricter bound on the number of derivatives needed for systems with multiple output equations. This method can be implemented in any symbolic algebra package with a Gröbner Basis function, and has the benefit of being simpler to implement than Ritt's pseudodivision algorithm, since differentiation and the Buchberger Algorithm are all that is required. The input-output equations of the system (1.1) can then be used for identifiability analysis, as demonstrated in [8, 9].